\input amstex
\input amsppt.sty
\def\Re{\operatorname{Re}}
\def\phi{\varphi}
\def\epsilon{\varepsilon}
\NoBlackBoxes

\topmatter
\title
Boundary behavior of the Kobayashi-Royden metric in smooth pseudoconvex domains
\endtitle

\author
Peter Pflug, W\l odzimierz Zwonek
\endauthor

\address
Carl von Ossietzky Universit\"at Oldenburg,
Institut f\"ur Mathe\-ma\-tik, Postfach 2503, D-26111 Oldenburg,
Germany
\endaddress

\email
pflug\@ mathematik.uni-oldenburg.de
\endemail

\address
Instytut Matematyki, Uniwersytet Jagiello\'nski, \L ojasiewicza 6, 30-348 Krak\'ow, Poland
\endaddress

\email
Wlodzimierz.Zwonek\@ im.uj.edu.pl
\endemail
\abstract
We show some lower estimates for the Kobayashi-Royden metric on a class of smooth bounded pseudoconvex domains.
\endabstract

\thanks Both authors were supported by
a DFG grant No 436 POL 113/106/0-2 (July 2009/September 2009) and by the the Research Grant No. N N201 361436 of the Polish Ministry of Science and Higher Education.
\endthanks

\thanks
\endthanks

\thanks
keywords: Kobayashi-Royden metric, pseudoconvex domain
\endthanks

\endtopmatter

\document

\subheading{1. Introduction and main results}

In this short note we discuss the problem of the boundary
behavior of the Kobayashi-Royden metric (mainly) in the normal
direction in smooth bounded pseudoconvex domains. We show two main
results. One of the results states in particular that the
Kobayashi-Royden metric in the normal direction in some class of
smooth bounded pseudoconvex domains is estimated from below by the
expression like $1/d_D^{7/8}(z)$ ($d_D(z)$ denotes the distance of
$z$ from the boundary of $D$) which improves a recent result of S.
Fu (see \cite{Fu~2009}) where the Author obtained the lower
estimate with the exponent $5/6$. On the other hand we show that a
careful study of a recent example of J. E. Fornaess and L. Lee
(see \cite{Lee~2009}) shows that the optimal exponent in the lower
estimate of the Kobayashi-Royden metric in the normal direction is
smaller than one (for $C^k$-smoothness, $k<\infty$) and we also
show some obstacles for the rate of the increase in the
$C^{\infty}$-case.

Recall that the Kobayashi-Royden metric has a localization property (see e.g. \cite{Roy~1971}, \cite{Gra~1975}); therefore, we lose no generality concentrating on domains defined globally. Recall also that one of reasons to study the boundary behavior of the Kobayashi-Royden metric is the problem of deciding whether any bounded smooth pseudoconvex domain is Kobayashi complete (see e.g. \cite{Jar-Pfl~1993}). The hope was that the Kobayashi-Royden metric in the normal direction explodes near the boundary as $1/d_D(z)$; it was one of the ideas to show that smooth bounded pseudoconvex domains are Kobayashi complete. However, after many years of uncertainty a recent example of Fornaess-Lee (\cite{For-Lee~2009}) showed that such a lower bound is not valid. More precisely, the example is the following.

\proclaim{Theorem 1 \rm{(see \cite{For-Lee~2009})}} For any given
increasing sequence $(a_{\nu})_{\nu}$, $a_\nu\to\infty$, of
positive numbers there is a bounded smooth pseudoconvex domain
$D\subset\Bbb C^3$ and a decreasing sequence
$(\delta_{\nu})_{\nu}$ with $\delta_{\nu}\to 0$ such that $$
\kappa_D(P_{\delta_{\nu}};n)\leq 1/(a_{\nu}\delta_{\nu}), $$ where
$P$ is a suitable point from $\partial D$,
$P_{\delta_{\nu}}=P-\delta_{\nu}n$, and $n$ is the unit outward
normal vector to $\partial D$ at $P$.
\endproclaim
In the other direction S. Fu showed that smooth bounded
pseudoconvex domains have the following lower estimate (our
formulation is weaker than the one in the original paper).

\proclaim{Theorem 2 \rm{(see \cite{Fu~2009})}} Let $D$ be a
bounded $C^3$-smooth pseudoconvex domain given by the formula
$D=\{r<0\}$, where $r$ is a $C^3$-smooth defining function (meaning
that its Levi form is semipositive definite on the complex tangent space
of any boundary point). Then there is a constant $c>0$ such that $$
\kappa_D(z;X)\geq c\frac{|\langle
\partial r(z),X\rangle|}{|r(z)|^{2/3}},\quad z\in D,\; X\in\Bbb C^n. $$
\endproclaim

Moreover, it follows from the paper of S. Fu that if we make
some additional assumption on the vector $X$ and points $z$ (for
instance $X$ is the unit outward normal vector to some
boundary point $P$ and $z$ lies on the line passing through $P$ in
the direction $X$) then in the above estimate we may replace the
exponent $2/3$ with $5/6$. We shall see in Theorem 4 that in many
cases the exponent $5/6$ may be replaced by $7/8$.

S. Fu also conjectured that in the class of smooth domains the lower estimate of the Kobayashi-Royden metric as in Theorem 2 may be taken of the form $1/d_D^{1-\epsilon}(z)$ with $\epsilon>0$ arbitrarily small. Note that the example of Fornaess-Lee shows that the exponent cannot be taken to be equal to one (equivalently $\epsilon$ cannot be equal to $0$).

However, the careful study of the example of Fornaess-Lee shows
that in the case of $C^k$-smooth domains an estimate as conjectured by S. Fu does
not hold. We have namely the following result.

\proclaim{Theorem 3} (1) For any positive integer $k$ there are a
$C^k$-smooth bounded pseudoconvex domain $D$ in $\Bbb C^3$, a
positive number $\epsilon$, and a decreasing sequence
$(\delta_{\nu})_{\nu}$ with $\delta_{\nu}\to 0$ such that $$
\kappa_D(P_{\delta_{\nu}};n)\leq 1/\delta_{\nu}^{1-\epsilon}, $$
where $P$ is a suitable point from $\partial D$,
$P_{\delta_{\nu}}=P-\delta_{\nu}n$, and $n$ is the unit outward
normal vector to $\partial D$ at $P$.

(2) For any $\alpha>0$ there are a $C^{\infty}$-smooth bounded
pseudoconvex domain $D$ in $\Bbb C^3$ and a decreasing sequence
$(\delta_{\nu})_{\nu}$ with $\delta_{\nu}\to 0$ such that $$
\kappa_D(P_{\delta_{\nu}};n)\leq \frac{1}{\delta_{\nu}(-\log
\delta_{\nu})^{\alpha}}, $$ where $P$ is a suitable point from
$\partial D$, $P_{\delta_{\nu}}=P-\delta_{\nu}n$, and $n$ is the
unit outward normal vector to $\partial D$ at $P$.
\endproclaim

Note that the above
theorem shows that even in the $C^{\infty}$ case the proof of the
Kobayashi completeness of the smooth bounded pseudoconvex domain
cannot go along the following lines: We prove that the
Kobayashi-Royden metric (in the normal direction) behaves like a
'regular' integrable function of $d_D(z)$. Therefore, in case that all smooth bounded pseudoconvex
domains are all Kobayashi
complete the proof would require a more subtle reasoning.

We can, however, say also something in the positive direction. Namely, we may slightly improve the estimate given in Theorem 2. Unfortunately, the
better estimate holds for smooth domains defined as sublevel sets of smooth plurisubharmonic defining function - recall that not all smooth bounded
pseudoconvex domains are locally sublevel sets of smooth plurisubharmonic functions (see \cite{For~1979}, \cite{Beh~1985}).

\proclaim{Theorem 4} Let $D=\{r<0\}$ be a bounded domain in $\Bbb
C^n$ where $r:U\mapsto\Bbb R$ is a $C^4$-smooth plurisubharmonic
defining function for $D$. Then there is a constant $C>0$ such
that $$ \kappa_D(z;X)\geq\frac{C|\langle
n(z),X\rangle|}{d_D^{7/8}(z)} $$ as $z$ tends to $\partial D$ and
the vectors $X$ are taken so that $||X||=o(1/d_D(z))|\langle
n(z),X\rangle|$, where $n(z)$ denotes the unit outward normal
vector to $\partial D$ at the point of $\partial D$ of the
smallest distance from $z$.
\endproclaim

Before we start the proofs recall the definition of the
Kobayashi-Royden (pseu\-do)\-met\-ric of a domain $D\subset\Bbb
C^n$ (for basic properties of the Kobayashi-Royden metric see
\cite{Jar-Pfl~1993}). $$ \kappa_D(z;X):=\inf\{\alpha>0:\text{
there is an $f\in\Cal O(\Bbb D,D)$ with } f(0)=z,\;\alpha
f^{\prime}(0)=X\}, $$ $z\in D,\; X\in\Bbb C^n$.

\subheading{2. Proofs}
We start with some preliminary considerations.

Let $D=\{r<0\}$ be a domain in $\Bbb C^n$, where $r:U\mapsto\Bbb R$
is a $C^{k+1}$-smooth plurisubharmonic defining function with
$r(0)=0$. Then (up to a linear isomorphism) the Taylor expansion
at $0$ of order $k$ of $r$ is of the following form $$
r(z)=\Re z_n+\sum_{j=2}^kQ_j(z)+R_k(z), $$ where
$Q_j(z)=\sum\sb{|\alpha|+|\beta|=j}a_{\alpha,\beta}^jz^{\alpha}\bar
z^{\beta}$ (note that then $a_{\alpha,\beta}^j=\bar
a_{\beta,\alpha}^j$). We may also write
$Q_j(z)=\tilde Q_j(z)+\hat Q_j(z)$,
where $$ \tilde
Q_j(z)=\sum\sb{|\alpha|+|\beta|=j,|\alpha|,|\beta|>0}
a_{\alpha,\beta}^jz^{\alpha}\bar z^{\beta},\quad \hat
Q_j(z)=2\Re\sum_{|\alpha|=j}a_{\alpha,\alpha}^jz^{\alpha}=:2\Re
H_j(z). $$

It follows from Taylor's formula that $\Cal L
R_k(z)|_{S^{2n-1}}=O(||z||^{k-1})$, where $\Cal L \tilde r(z)(X)$
is the Levi form of $\tilde r$ at the point $z$ in direction of
$X$; $S^{2n-1}$ means the $2n-1$-dimensional sphere. Consequently,
$\Cal L R_k(z)(z)=O(||z||^{k+1})$.

An easy calculation gives the following formula $$
\multline \Cal L
Q_j(z)(z)=\sum_{\nu=0}^j\nu(j-\nu)\sum_{|\alpha|=\nu,|\beta|=j-\nu}a_{\alpha,\beta}^jz^{\alpha}\bar
z^{\beta}=\\
\sum_{\nu=1}^{j-1}\nu(j-\nu)\sum_{|\alpha|=\nu,|\beta|=j-\nu}a_{\alpha,\beta}^jz^{\alpha}\bar
z^{\beta}.
\endmultline
$$
In particular, $\Cal L Q_2(z)(z)=\tilde Q_2(z)$, $\Cal L Q_3(z)(z)=2\tilde Q_3(z)$.

In the sequel we shall denote by $C_j$ different constants that depend only on the domain $D$.

\demo{Proof of Theorem 4} We leave the notation as above and we make use of the above considerations.

First we prove the desired estimate but with the exponent equal to $3/4$ (instead of $7/8$).

The assumptions of the theorem imply that for $z\in U$ we have the following estimate
$$
\tilde Q_2(z)+ 2\tilde Q_3(z)+C_1||z||^4\geq 0,
$$
which together with the property
$$
\min\{\tilde Q_2(z)+\tilde Q_3(z),\tilde Q_2(-z)+\tilde Q_3(-z)\}\geq\min\{\tilde Q_2(z)+2\tilde Q_3(z),\tilde Q_2(-z)+2\tilde Q_3(-z)\}
$$
gives for $z$ close to $0$ the inequality
$$
r(z)\geq\Re z_n+Q_2(z)+Q_3(z)-C_2||z||^4\geq\Re z_n+2\Re(H_2(z)+H_3(z))-C_3||z||^4.
$$
Therefore, shrinking $U$ if necessary, we have the following inclusion
$$
D\subset\{\Re z_n+2\Re(H_2(z)+H_3(z))-C_3||z||^4<0\}.
$$

Now for $\delta>0$ small enough and $X\in\Bbb C^n$, $X\neq 0$,
take $\phi\in\Cal O(\Bbb D,D)$ such that $\phi(0)=(0,\ldots,0,-\delta)$,
$\kappa\phi^{\prime}(0)=X$ where $\kappa>0$.

Note that $||\phi(\lambda)-\phi(0)||\leq C_4|\lambda|$.
For $r\in(0,1)$  define $\psi_r(\lambda):=\phi(r\lambda)$. Then $||\psi_r(\lambda)||\leq\delta+C_4 r$, $\lambda\in\Bbb D$.

Define $\Psi(z):=z_n+2H_2(z)+2H_3(z)$. Put $\phi_r:=\Psi\circ\psi_r$. Then $\phi_r(\Bbb D)\subset\{\lambda\in\Bbb C:\Re\lambda<C_5 (\delta+r)^4\}=:S_{r,\delta}$
for $\delta$ small enough.

Therefore,
$$
\frac{r}{\kappa}|X_n(1+O(\delta))+\sum_{j=1}^{n-1}O(\delta)X_j|=\frac{r}{\kappa}\left|\Psi^{\prime}(0,\ldots,0,-\delta)X\right|=
|\phi_r^{\prime}(0)|\leq C_6(\delta+(\delta+r)^4),
$$
where the last inequality follows easily from the formula for the Kobayashi-Royden metric for $S_{r,\delta}$.

Substitute $r=\delta^{1/4}$ for $\delta$ small enough.
Then we get the following lower estimate
$$
\kappa\geq C_7 \delta^{1/4}\frac{|X_n|(1+\alpha(\delta))}{\delta+(\delta+\delta^{1/4})^4},
$$
where $\alpha(\delta)\to 0$ as $\delta\to 0$ (here we use the fact that
we choose $X$ such that $||X||=o(1/\delta)|X_n|$). Consequently, we get the following lower estimate
$$
\kappa_D((0,\ldots,0,-\delta);X)\geq C_6|X_n|/\delta^{3/4}
$$
as $\delta$ tends $0$ and the vectors are taken so that $o(1/\delta)|X_n|\geq||X||$.

Our aim is to show that we may replace the exponent $3/4$ with
$7/8$. Keeping in mind the above estimate and leaving the same
notation as above we have the following inequality
$$
||\phi(\lambda)-\phi(0)-\frac{\lambda}{\kappa}X||\leq
C_7|\lambda|^2. $$

Proceeding as before we get for $|\lambda|\leq r$, $\delta$ small
enough $$
||\phi(\lambda)||\leq\delta+\frac{rC_8}{\kappa}||X||+C_7r^2 $$ or
$||\psi_r(\lambda)||\leq \delta+\frac{rC_8}{\kappa}||X||+C_7r^2$,
$\lambda\in\Bbb D$.

Since without loss of generality we may assume that $||X||$ is bounded from above (or even equal to one)
proceeding exactly as before we get the following inequality
$$
\kappa\geq\frac{C_9r|X_n|(1+\alpha(\delta))}{\delta+(\delta+\frac{r}{\kappa}+r^2)^4}.
$$
Since we already know that $\kappa\geq C_{10}/(\delta^{3/4})$ (at this place we need the first part of the proof),
putting $r=\delta^{1/8}$ we get the following estimate
$$
\kappa_D((0,\ldots,0,-\delta);X)\geq\frac{C_{11}|X_n|}{\delta^{7/8}}
$$
as $\delta$ tends to $0$ and $X$ satisfies the inequality $||X||\leq o(1/\delta)|X_n|$, which finishes the proof of the theorem.
\qed
\enddemo

\subheading{Remark 5} Consider $r$ to be defined near $0$ as follows
$r(z):=\Re z_2+p(z_1)$ where
$$
\multline
p(z_1):=2m^2z_1^{m+l}\overline z_1^{m-l}+4(m^2-l^2)|z_1|^{2m}+2m^2z_1^{m-l}\overline
z_1^{m+l}=\\
|z_1|^{2m}(\Re{4m^2e^{i2l\theta}+4(m^2-l^2)}),
\endmultline
$$
$m/2\leq l<m$ (here $z_1=e^{i\theta}|z_1|$). Since $p$ is a subharmonic function such that for some values of $\theta$ the last factor in the
formula is negative (the example is taken from \cite{Las~1988}) we see that we cannot hope to repeat the reasoning from the proof of
Theorem 4 for general $k$ - even the case $k=2m=4$ ($m=2$, $l=1$) encounters an obstacle. In other words the above method of the proof does not
give a better lower estimate. Nevertheless, we think that the lower estimate in the normal direction of the Kobayashi-Royden metric near the boundary
of a $C^{k+2}$-smooth pseudoconvex domain
may be of the form $1/d_D^{1-1/(2k)}(z)$, which would mean that in the case of infinitely smooth bounded pseudoconvex
domain the estimate with the exponent arbitrarily close to $1$ (as suggested by S. Fu) may hold.

\subheading{Remark 6} In the proof of Theorem 4 one needs twice
the same reasoning. However, instead of repeating it twice one may
use a result of S. Fu (to get the lower estimate of $\kappa$ of
the form $1/\delta^{2/3}$ - in fact it is sufficient to have the
estimate of the form $1/\delta^{1/8}$). However, at the present
form the proof is more self contained. Therefore, the authors decided to leave it in the present
form.

\demo{Proof of Theorem 3} As mentioned earlier the domain which satisfies the properties claimed in the theorem was constructed in
\cite{For-Lee~2009}. Therefore, we recall the construction from there (keeping the notation from there, too). To get the proof of the theorem we have to add some estimates (mostly for derivatives) of the defined functions and also at some places, for simplicity of calculations, we make some special choice of some sequences.

First let us make some comments. Note that the procedure works not only
for $r_{n+1}=\frac{r_n^2}{a_n}$ but also under the assumption that $r_{n+1}\leq\frac{r_n^2}{a_n}$, also the choice of $A_k$ may be done with the equality replaced by the inequality. Consequently, the series that we shall choose can be replaced by any subseries.
%In other words the condition that may be verified is that the appropriate sequence must be convergent to $0$.

Below we shall write some inequalities for norms. Such an inequality: $l_n\leq m_n$ means that $\limsup\frac{l_n}{m_n}<\infty$. The meaning of the equality is analogous (inequalities in both directions hold). The norms of functions are meant to be the supremum norms of functions (on some sets).

At first stage we repeat the definition of a sequence of
subharmonic functions which is then adopted to the construction of
a sequence of subharmonic functions of two variables defining a
three-dimensional example. As mentioned earlier the construction
follows entirely from \cite{For-Lee~2009}.

At first we assume the existence of sequences $(a_n)_n$, $(r_n)_n$ such that the sequence $(a_n)_n$ is increasing to infinity and $r_{n+1}\leq \frac{r_n^2}{a_n}$. We shall fix the sequences later.

We define $$ u_n(z):=1/8-\Re z+\frac{\log|z|}{4\log a_n},\;
z\in\Bbb C. $$ Then we $$ R_n(z):= \cases
\max\{u_n(z),0\},&\Re z\leq b_n\\ u_n(z),& \Re z>b_n
\endcases,
$$
where $0< b_n\leq 1$ is the smallest positive number such that
$$
1/8-b_n+\frac{\log b_n}{4\log a_n}=0.
$$
At first we are interested in the norm of $R_n$ on a closed disc of radius $Ma_n/r_n$ (for some fixed $M>1$). It is estimated from above by $a_n/r_n$.

We define
$$
\tilde R_n(z):=\int_{\Bbb C}R_n(z-\epsilon_n w)\chi(w)d\mu(w)
$$
for some $0<\epsilon_n<r_n/2$, where $\mu=dxdy/m$. Here $m=\int_{\Bbb C} \chi(z)dxdy$ and $\chi:\Bbb C\mapsto[0,\infty)$ is a non-constant $C^{\infty}$ radial function such that $0\leq\chi\leq 1$ and $\chi(z)=0$ for $|z|\geq 1$.

Then $||\tilde R_n^{(k)}||_{B(0,M)}\leq \left(1/r_n\right)^k \cdot a_n/r_n$. Then we put $\rho_n(z):=\tilde R_n(\frac{a_nz}{r_n})$.

We therefore have
$$||\rho_n^{(k)}||_{B(0,M)}\leq \left(\frac{a_n}{r_n^2}\right)^k\frac{a_n}{r_n}$$.

At this place we fix the sequences. We put $r_n=1/a_n$. We also want
to have $r_{n+1}:=\frac{r_n^2}{a_n}=r_n^3$. In other words we may
choose $r_n:=r_1^{3^n}$.
Now fix for a while $\epsilon\in(0,1/3)$. And put
$a_n=\delta_n^{-\epsilon}$. So finally, the choice of the numbers
is the following $a_n=a^{\epsilon 3^n}$, where $a>1$ is fixed,
$r_n=(1/a)^{\epsilon 3^n}$, $\delta_n=(1/a)^{3^n}$.
But the construction  needs also additional number
$A_n=1/2+a_n/r_n+\log(1/r_n)/4\log a_n$ such that
$$
\delta_n\leq\frac{\delta_{n-1}}{A_n\frac{1}{2^n}}.
$$

With our choice of numbers we get that
$A_n=a^{2\epsilon 3^n}$ (in the asymptotic sense). The construction needs also that
$\delta_n\leq\frac{\delta_{n-1}}{A_n 2^n}$. Since the inequality
must hold asymptotically (it follows from the reasoning) it is
sufficient to see that for large $n$
$$
a^{-3^n}\leq\frac{1}{a^{3^{n-1}+2\epsilon 3^n}2^n},
$$
which holds for $\epsilon$ as above.

At this place the one-dimensional function $\rho$ is defined as
$$
\rho(z):=\sum\sb{n=1}\sp{\infty}\delta_n\rho_n(z),
$$
which defines a $C^{k}$-smooth function under the assumption
$$ \sum\sb{n}\delta_n||\rho_n^{(k)}||_{B(0,M)}<\infty. $$
In other
words this gives the condition
$\sum\sb{n}\delta_n\frac{a_n^{k+1}}{r_n^{2k+1}}<\infty$. But the
last series is $\sum a^{((3k+2)\epsilon-1)3^n}$, which is finite
when $\epsilon<1/(3k+2)$.

Now we move to the construction of the proper function $\tilde\rho$.

We define $V:=\{(s,t)\in\Bbb C^2:s^2-t^3=0\}$. We want to have $\tilde\rho_n(s,t)=\rho_n(s/t)=\rho_n(\zeta)$ if $(s,t)=(\zeta^3,\zeta^2)\in V$.

Let $\tilde r_n:=r_{n+1}^3$ and put $B_n:=B(0,\tilde
r_n)\subset\Bbb C^2$, we also put $B_n^{\prime}:=B(0,3/4\tilde
r_n)$. Then one may choose a small neighborhood $U_n$ of $V$ such
that the projection $\pi:U_n\mapsto V$ is well-defined on
$U_n\setminus B_n^{\prime}$ (the formula is the following
$\pi(s,t):=(s,s^{2/3})$ with a properly chosen branch of the
power). We put $U_n:=\{p\in\Bbb C^2:||p-\pi(p)||<d_n^2\}$, where
one may choose $d_n=r_{n+1}^3=r_n^9$ (asymptotically in the above
mentioned sense).

Then $||\pi^{(k)}||_{U_n\setminus B_n^{\prime}}=r_{n+1}^{-k}=\frac{1}{r_n^{3k}}$.

We define $\tilde\rho_n:=\rho_n\circ\pi$ on $U_n\setminus B_n$ and
we may extend $\tilde\rho_n$ to a $C^{\infty}$-smooth on $B_n\cup
U_n$ letting it be equal $0$ on $B_n$.

Now we note the next estimate 
$$
||\tilde\rho_n^{(k)}||_{(U_n\setminus B_n)\cap B(0,M)}\leq\frac{1}{r_n^{6k+2}}.
$$

Let $\chi:\Bbb R\mapsto[0,1]$ be a $C^{\infty}$-smooth function, equal to $1$ on $[0,1/2]$ and equal to $0$ on $[1,\infty)$.

Then we define another smooth extension of $\tilde\rho_n$ on $\Bbb C^2$ by the formula
$$
p_n(z):=\cases
0,& z\in B_n,\\
\tilde\rho_n(z),&z\in U_n\setminus B_n,||z-\pi(z)||\leq\frac{d_n^2}{2}\\
\tilde\rho_n(z)\chi\left(\frac{||z-\pi(z)||^2}{d_n^2}\right),&z\in U_n\setminus B_n,\frac{d_n^2}{2}\leq||z-\pi(z)||^2\leq d_n^2,\\
0,&z\not\in U_n\cup B_n
\endcases
$$
Then one may verify that $||p_n^{(k)}||_{B(0,M)}\leq\frac{1}{r_n^{9k+2}}$. Now we take $C_n\geq 0$ such that
$\Cal L p_n(z)(X)\geq-C_n||X||^2$. It follows that we may take (asymptotically)
$$
C_n=\frac{1}{r_n^{20}}
$$
-- use the estimate for the norm of $p_n^{\prime\prime}$.

Put $$ A_n:=\{z\in U_n\setminus
B_n:\frac{d_n^2}{2}\leq||z-\pi(z)||^2\leq d_n^2\} $$  and  $$
q(s,t):=e^{||(s,t)||^2}|s^2-t^3|^2,\;(s,t)\in\Bbb C^2. $$ Note
that $\Cal Lq(s,t)(X)\geq |s^2-t^3|^2||X||^2$. Therefore, if we
take (asymptotically) $c_n=d_n^2=r_n^{18}$, then $$ \Cal
Lq(z)(X)\geq c_n||X||^2. $$

Put $K_n=\frac{1}{r_n^{38}}$ (asymptotically). Then $-C_n+K_nc_n\geq
0$. Consequently, $\tilde r_n=p_n+K_nq$ is plurisubharmonic on $\Bbb
C^2$. Certainly, $||\tilde r_n^{(k)}||_{B(0,M)}\leq
\max\{\frac{1}{r_n^{9k+2}},\frac{1}{r_n^{38}}\}=:\frac{1}{r_n^{m_k}}$
(note that $m_k=38$, $k=1,\ldots,4$ and $m_k=9k+2$, $k\geq 4$).
Now the condition on $C^{k}$-smoothness of the example from
\cite{For-Lee~2009} follows from the $C^k$-smoothness of
$$
\tilde\rho:=\sum_n\delta_n \tilde r_n,
$$
which is satisfied if
$$
\infty>\sum_n\delta_n\frac{1}{r_n^{m_k}}=\sum_na^{m_k\epsilon-1}.\tag{$\ast$}
$$
The last inequality completes the proof with arbitrary $\epsilon\in(0,1/m_k)$.

To complete the construction recall that Fornaess and Lee defined the domain as follows
$$
D:=\{(s,t,w)\in\Bbb C^3:\Re w+\tilde\rho(s,t)<0\}\cap B(0,2).
$$

Let us now move to the second part of the theorem.

We leave all the relations between the
numbers $a_n,\delta_n,r_n$ with one exception. Namely, put
$a_n=(-\log\delta_n)^{\alpha}$, $\alpha>0$. Explicitly we have
$\delta_n=(1/a)^{3^n}$, $r_n=1/a_n$. Then the convergence of the final sequence \thetag{$\ast$} (with just
introduced $\delta_n$ and $r_n$) is easily satisfied. And although the relation
$r_{n+1}\leq\frac{r_n^2}{a_n}$ is not satisfied now, it is easy to see that instead of taking the whole sequence
while defining the function $\tilde\rho$ we may also choose an arbitrary subsequence which easily guarantees that the
desired inequality is satisfied. One may easily prove that choosing these relations the number $A_n$ satisfies the desired inequality as well.

Above considerations lead us to the following relation being sufficient
for the construction of a $C^{\infty}$-smooth domain with the
boundary behaviour of the Kobayashi-Royden metric in the normal
direction equal to $1/(\delta_n(-\log\delta_n)^{\alpha})$:
$$
\sum\delta_na_n^{k}<\infty  \;\text{ for any positive integer $k$}.
$$
The last condition is, as one may verify, satisfied.
\qed
\enddemo

\enddocument